\documentclass[journal]{IEEEtran}
\usepackage{cite}
\usepackage{textcomp}
\usepackage{tabulary}
 \usepackage{blindtext}
\usepackage{scrextend}
\usepackage{float}
\usepackage[usenames, dvipsnames]{color}
\usepackage{enumerate}
\usepackage{subfigure}
\ifCLASSINFOpdf
   \usepackage[pdftex]{graphicx}
\else
 \usepackage[dvips]{graphicx}
\fi
\usepackage[cmex10]{amsmath}
\usepackage{url}
\begin{document}
\title{A Protection Method in Active Distribution \\ Grids with High Penetration of Renewable \\ Energy Sources}

\author{J.K.~Wang,~\IEEEmembership{Member,~IEEE,} and Christian Moya Calderon,~\IEEEmembership{Student Member,~IEEE}
\thanks{J.K.~Wang and C.~Moya Calderon are with the Department
of Electrical and Computer Engineering, The Ohio State University, Columbus,
OH, 43210 USA e-mail: \{wang.6536, moyacalderon.1\}@osu.edu.}}

\markboth{}%
{Shell \MakeLowercase{\textit{et al.}}: Bare Demo of IEEEtran.cls for Journals}

\maketitle
\begin{abstract}
A protection method in active distribution networks is proposed in this paper. In active distribution systems, fault currents flow in multiple directions and presents a varying range of value, which poses a great challenge of maintaining coordination among protective devices on feeders. The proposed protection method addresses this challenge by simultaneously adjusting DG's output power and protection devices' settings in pre-fault networks. Comparing to previous protection solutions, the proposed method considers the influences from renewable DG's intermittency, and explores the economic and protection benefits of DG's active participation. The formulation of proposed method is decomposed into two optimization sub-problems, coupling through the constraint on fuse-recloser coordination. This decomposed mathematical structure effectively extinguishes the non-linearity arising from reclosers' time-current inverse characteristics, and greatly reduces computation efforts. 
\end{abstract}
\vspace{-4mm}
\section{Nomenclature}
\noindent {\bf{\it Variables}}
\begin{labeling}{alligator}
\item [$I^f_P, I^f_B$]  fault current seen by the primary and backup protective device,
\item [$\Delta I^f_{FR,j}$] fault current disparity between recloser at node $j$ and its lateral fuse,
\item [${\Delta} I^f_{RR,j}$] fault current disparity between the primary recloser at node $j$ and its backup recloser $j-1$,
\item [$I^f_{max}, I^f_{min}$]  maximum and minimum fault current,
\item [$I^f_{R,j}$] fault current seen by the recolser at node $j$,
\item [$I^f_{F,i}$] fault current seen by the fuse at lateral $i$,
\item [$I^f_{G,i}$] fault current contributed by DG connected at lateral $i$,
\item [$P_i,Q_i$] real and reactive power flowing on the feeder section between lateral $i$ and lateral $i+1$ ,
\item [$P_{D,i},Q_{D,i}$] total real and reactive power demand of all loads on lateral $i$,
\item [$P_{G,i},Q_{G,i}$] real and reactive power output from DG at lateral $i$,
\item [$V_i$] terminal voltage magnitude of lateral $i$,
\end{labeling}
\vspace{2mm}
\noindent {\bf{\it Functions}}
\begin{labeling}{alligator}
\item[$T_{R,j}$] Time-Current Inverse (TCI) curve of recloser $j$,
\item[$T_P, T_B$] response time of the primary and backup protective device, which is calculated by the TCI curve $T_P(I^f_P)$ and $T_B(I^f_B)$ .
\end{labeling}

\vspace{2mm}
\noindent {\bf{\it Parameters}}
\begin{labeling}{alligator}
\item[$r_i,x_i$]  resistance and reactance of the feeder section between lateral $i$ and $i+1$,
\item[${\Delta}T_{FR}$] coordination margin between the recloser and its lateral fuse,
\item[${\Delta}T_{RR}$]  coordination margin between two reclosers.
\item[$D_{R,j}$] time dial setting of recloser $j$,
\item[$I^p_{R,j}$] pickup current of recloser $j$,
\end{labeling}

\vspace{2mm}
\noindent {\bf{\it Sets}}
\begin{labeling}{alligator}
\item[$\mathcal{J}$] bus/nodes on the main feeder,
\item[$\mathcal{I}$] laterals and DG inter-ties tapped on the main feeder,
\item[$\mathcal{I}^{j}$] laterals tapped on the feeder section between node $j$ and $j+1$,
\item[$\mathcal{I}^u_{j},\mathcal{I}^d_{j}$] laterals upstream and downstream to node $j$,   
\end{labeling}

\section{Introduction}
\begin{figure}[!t]
\centering
	\includegraphics[width=0.48\textwidth]{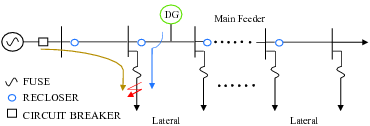}
\caption{Protection configuration of a distribution network.}
\label{fig: Fig1}
\vspace{-4mm}
\end{figure}

\IEEEPARstart{I}{ncreasing} penetration of Distributed Generation (DG) at distribution level of utility grids (generally radial up to $35kV$) poses challenge of maintaining coordination in legacy protection systems. As illustrated in Fig.~\ref{fig: Fig1}, conventional distribution grid is characterized by a single source (substation) feeding a network of downstream feeders. Legacy protection systems are designed assuming the fault current flowing in one direction. After connecting DG, the system becomes multi-source fed, and protective devices' coordination may not hold properly \cite{Girgis,Brahma,Chaitusaney}. In addition, fault current values in an active distribution grid can vary over a wide range due to DG's intermittent output power \cite{ieee2004impact, de2012impact}. This makes setting protective devices difficult: settings that are conservative may sacrifice protection sensitivity, and settings that are strict may risk protection security and lead to mal-operation \cite{ieee2004impact}.

To retrieve protection coordination, many methods have been proposed and can be classified as (i) limiting DG's fault current contribution \cite{tang2005application, El-Khattam}; (ii) deploying microprocessor-based recloser/relay \cite{Chaitusaney,Brahma,Abdelaziz}; (iii) deploying advanced protection schemes that are not overcurrent \cite{dewadasa2011protection, haron2012review, sortomme2010microgrid}. An advantage of methods in (i) and (iii) is that they remove the need of coordination \cite{ieee2004impact, Barker}; with the support from advanced communication and metering, some method can even achieve ``set-less'' performance \cite{sortomme2010microgrid}. However, these methods either require high capital investment, or conflict with the political goal and environmental requirement by limiting DG's penetration level. Methods in (iii) are more suitable to compact electrical systems, such as a micro-grid, in which loads are aggregated under buses. For wide-spread distribution networks, in which loads are tapped along feeders, closed-zone protection schemes (i.e. differential protection) are difficult to achieve full coverage.

Methods in (ii) provide promising solutions. They utilize the flexibility provided by microprocessor-based recloser/relays, and choose appropriate settings by considering DG's installation location and types \cite{Brahma, Abdelaziz, chaitusaney2005impact}. In case DG units are disconnected temporarily, settings of recloser/relays can be adapted to the system's configuration in a timely manner \cite{Abdelaziz}.

One problem that has not been addressed by the existing studies in (ii) is the limitation of recloser's adaptivity. With more DG installed, solely adjusting recolsers' settings would not be sufficient to maintain protection coordination. To address this problem, we propose a protection method that simultaneously adjusts the settings of recloser/relays and DG's output power in pre-fault distribution networks. The proposed method ensures secure operation of protective devices and fast fault clearing in active distribution grids. System operation benefits are increased by maximizing DGs' total power output. Higher DG penetration capacity can be achieved by smoothening the fault current fluctuation from renewable DG.  

The rest of the paper is organized as following: Section III reviews protection fundamentals and DG's impact on protection in active distribution grids. Section IV discuss DG's fault current contribution by examining different DG types and considering DG's connection requirements. Section V presents the proposed protection method, its mathematical formulation and algorithm. In Section V, the proposed method's benefits are verified on the IEEE 37-bus system. The contributions of this paper are summarized in Section VI.

\section{Background}
Conventional distribution grids are radial in nature, so that faults can be isolated conveniently from its upstream networks. Fig.~\ref{fig: Fig1} shows a typical protection configuration. Before DG's integration, fault current is fed from a single source and decreases along the feeder, flowing to the fault location \cite{Barker}.  
The main feeders are protected by reclosers and relays. A distribution network usually installs a few reclosers along the feeders, and one relay at the head of each feeders. Along a main feeder, laterals are tapped to supply loads. Laterals are protected by fuses. With distribution system becomes more automated, some laterals are installed with reclosers \cite{Fumilayo}. However, the large number of laterals and high cost of lateral reclosers make the fuses, which have much lower cost and much less application complicacy, the top choice for lateral protection. 
%
\vspace{-2mm}
\subsection{Protection coordination requirements}
Protection coordination directly affects system reliability. For feeders' protection, coordination is required among (i) fuse-fuse; (ii) fuse-recloser, and (iii) recloser-recloser. The last category includes recloser-relay. During the restoration post-fault, coordination may be required among protective devices of feeders , loads (primarily induction motors) and DG's reclosers \cite{anderson, Kersting}. In this paper, we restrict the discussion to the pre-fault and fault clearing process.

Two concepts defining coordination performance are {\it coordination range} and {\it coordination margin}. A pair of protective devices must be coordinated for all the fault currents within $I^f_{min}$ and $I^f_{max}$. This range of current is called the coordination range \cite{Girgis}. A time margin must be preserved for the response time between the primary and backup protective devices . This time is called coordination margin \cite{anderson}. Mathematically, the primary and backup devices are coordinated if and only if,

\begin{align} \label{cod_range}
\begin{split}
  T_P(I^f_{min})\leq T_B(I^f_{min})\\
       T_P(I^f_{max})\leq T_B(I^f_{max})
\end{split}
\end{align}  
and 
\begin{equation} \label{cod_margin}
T_B(I^f_B) - T_P(I^f_P)  \geq \Delta T; \forall I^f_B, I^f_P\in [I^f_{min},I^f_{max}]
\end{equation} 

\begin{figure}[t!]
\vspace{-4mm}
\setlength{\abovecaptionskip}{-3 pt}
\hspace*{-6mm}
\includegraphics[width=0.53\textwidth]{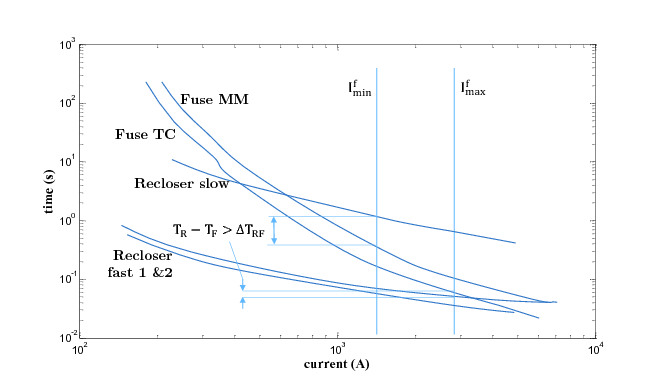}
\centering
\caption{Coordination margin and coordination range for a recloser and its lateral fuse. The reclosing sequence is fast-fast-slow.}
\label{fig: fr1}
\vspace{-4mm}
\end{figure}

\begin{figure}[t!]
\vspace{-4mm}
\setlength{\abovecaptionskip}{-10 pt}
\hspace*{-6mm}
\includegraphics[width=0.52\textwidth]{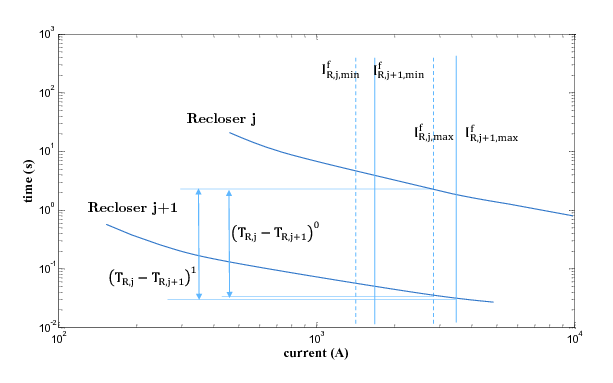}
\centering
\caption{Backup delay increasing due to DG's installation. The backup response delay is $(T_{R,j}-T_{R,(j+1)})^0$ before DG's installation, and increases to $(T_{R,j}-T_{R,(j+1)})^1$ after DG units connected to the feeder section between node $j$ and $j+1$. The two reclosers see different fault currents after DG's installation.}
\label{fig: rr1}
\vspace{-4mm}
\end{figure}

Fig.~\ref{fig: fr1} and Fig.~\ref{fig: rr1} illustrates the coordination range for fuse-recloser and recloser-recloser. In a passive distribution system (i.e. before DG's integration), the primary and backup devices see the same fault current, $I^f_P = I^f_B$.  

The coordination for fuse-fuse and recloser-recloser follows the principle that upstream devices backup downstream devices. In distribution grids, $70\%-95\%$ faults are temporary faults and can be cleared by reclosing operation \cite{glover2011power, Kersting, ieee2004impact}. Relays are used to clear up the permanent faults that cannot be cleared by downstream reclosers and backup for all protective devices on feeders.

The recloser-fuse coordination scheme varies by utilities' practice. Basically, two coordination schemes are deployed, fuse saving and fuse sacrificing  \cite{ieee2004impact}. In fuse sacrificing scheme, a fuse picks up a lateral fault before its upstream recloser; whereas in fuse saving scheme, a recloser first attempts to clear the fault by fast tripping. If the fault still persists, it is a permanent fault and is cleared up by the lateral fuse. In case the fault is located between the fuse and recloser, the recloser takes another slow trip after a fuse blown. Most recloser-fuse pairs are coordinated under fuse saving scheme. As shown in Fig.~\ref{fig: fr1},  a recloser's fast tripping curve is  coordinated with a fuse's Minimum Melting (MM) curve, and its slow tripping curve is coordinated with a fuse's Total Clearing (TC) time. In practice, a reclosing sequence takes various form, for example, one fast two slow (F-F-S) or two fast two slow (F-F-S-S) \cite{anderson, Chaitusaney, Girgis, Fazanehrafat}. 

\vspace{-2mm} 
\subsection{DG's impact on coordination and solutions}
In an active distribution grid, fault currents could flow from multiple directions. Hence, reclosers and relays must be retrofit with directional schemes to retrieve coordination \cite{Barker, Girgis}. Since fuses cannot sense currents' direction, if a large amount of DG are connected to laterals, fuse-fuse coordination cannot be retrieved. In this case, fuse needs to be replaced with more advanced protective devices \cite{Fumilayo}. In this paper, we do not attempt to retrieve fuse-fuse coordination and assume DG's are mainly integrated on feeders.    

\subsubsection{Fuse-recloser coordination}
In addition to fault current direction, connecting DG on feeders may disrupt fuse-recloser coordination in two ways: (i) increasing the maximum fault current on feeders; (ii) causing a fault current disparity seen by a coordinated pair \cite{Girgis, Brahma}. In the first case, DG's fault current contribution makes the fault current seen by the fuse exceeding the coordination range. The fault current seen by the recloser and fuse are the same, $I^f_{R,j}=I^f_{F,j}$. Coordination condition \eqref{cod_range} fails if \[I^f_{R,j}+\Delta I^f_{FR,j} > I^f_{max}.\]

In the second case, DG's fault current contribution makes the fuse response time falling out of the coordination margin. Coordination condition \eqref{cod_margin} fails if \[T_{R,j}(I^f_{max}) < T_{F,i}(I^f_{max} + \Delta I^f_{FR,j}).\]  

%

In either case, $I^f_{max}$ for the fuse is increased due to DG's fault current contribution. This fault current increment disrupts the coordination between a recloser' fast tripping curves and a fuse's MM curve. Certainly, DG's fault current contribution also lead to increased $I^f_{min}$. However, this change still will not disrupt the coordination between a recloser's slow tripping curves and a fuse's TC curve. 
\subsubsection{Recloser-recloser coordination}
DG's penetration is unlikely to disrupt recloser-recloser coordination \cite{Abdelaziz}. However, DG's fault current contribution could delay the response from the backup recloser.  Fig.~\ref{fig: rr1} illustrates this effect. For a pair of reclosers with DG units connected to the feeder section in between, the primary recloser sees more fault current, $I^f_{R,j}=I^f_{R,j-1}+\Delta I^f_{RR,j}$, and responses faster. The backup time delay is increased by \[T_{R,j-1}(I^f_{R,j-1}+\Delta I^f_{RR,j})-T_{R,j-1}(I^f_{R,j-1}).\] (The coordination margin is enlarged by $T_{R,j-1}(I^f_{max}+\Delta I^f_{RR,j})-T_{R,j-1}(I^f_{max})$. Hence, coordination condition \eqref{cod_margin} always holds.) With more DG penetrated, this backup delay could become large. Late backup actions will cause damage to equipment. 

\begin{figure}[t!]
\vspace{-3mm}
\setlength{\abovecaptionskip}{-10 pt}
\hspace*{-6mm}
\includegraphics[width=0.53\textwidth]{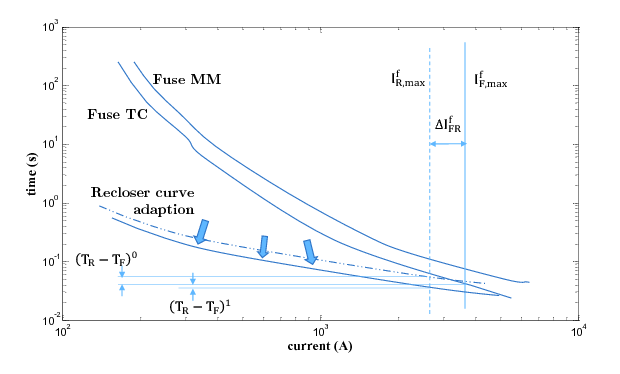}
\centering
\caption{Retrieving fuse-recloser coordination by adapting recloser's fast curve. After DG's installation, the fuse picks up the fault current before the recloser. It is observed at maximum fault current, $(T_R-T_F)^0>0$. By adapting the recloser fast tripping curve to DG's fault current contribution ${\Delta} I^f_{FR}$, the recloser picks up the fault before the fuse, $(T_R-T_F)^0<0$. }
\label{fig: fr3}
\vspace{-4mm}
\end{figure} 

\subsubsection{Deficiencies of existing adaptive protection methods}
The solution to the above problems depends on two factors: DG's fault current contribution $I^f_G$ and protective devices' TCI curves. 

A few solutions by adjusting the recloser/relay's TCI curves are proposed in previous studies \cite{Abdelaziz, Chaitusaney, Brahma}. According to the anti-islanding requirements, DG must be disconnected on detection of abnormalities in system operation \cite{basso2004ieee}. The required time for DG to be disconnected varies from a few cycles to 5 seconds \cite{kumpulainen2004analysis, de2012impact}. Based on the current anti-islanding practice, it is common to assume that DG's are disconnected after the first trip of reclosers \cite{Brahma, Chaitusaney, Barker}.

Therefore, to retrieve fuse-recloser coordination, we only need to coordinate a recloser's first tripping curve and a fuse's MM curve. The recloser's first fast curve should be moved downward until the coordination range and coordination margin are restored, as shown in Fig.~\ref{fig: fr3}.

\begin{figure}[t!]
\vspace{-3mm}
\setlength{\abovecaptionskip}{-3 pt}
\hspace*{-3mm}
\includegraphics[width=0.52\textwidth]{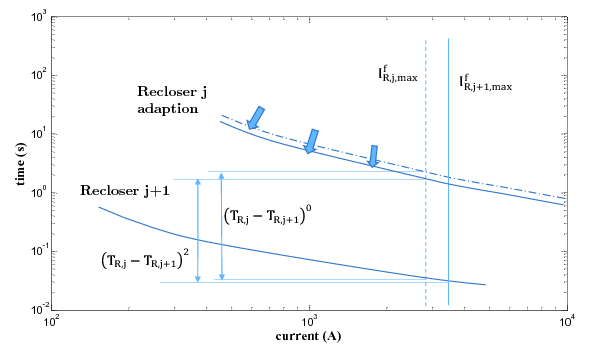}
\centering
\caption{Reducing backup delays between reclosers. By adapting Recloser $j$. The backup delay is restored to $(T_{R,j}-T_{R,j+1})^2$, which is the same as the backup delay before DG's installation$(T_{R,j}-T_{R,j+1})^0$.}
\label{fig: rr2}
\vspace{-1mm}
\end{figure}

To reduce extra backup delay among reclosers, the TCI curves of the backup recloser is moved downward until the coordination margin is reached, as shown in Fig.~\ref{fig: rr2}. In radial distribution systems, an upstream recloser backups all its downstream devices \cite{Abdelaziz, Brahma}.  

\begin{figure}[t!]
\vspace{-3mm}
\setlength{\abovecaptionskip}{-3 pt}
\hspace*{-3mm}
\includegraphics[width=0.49\textwidth]{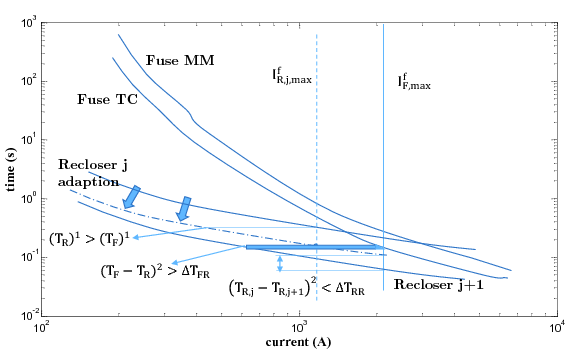}
\centering
\caption{Limitation of existing solutions to retrieving protection coordination. Coordination is lost after DG's installation, $(T_R)^1> (T_F)^1$. By adapting recloser $j$, coordination can be retrieved between recloser $j$ and the fuse, $(T_F-T_R)^2 < {\Delta}T_{FR}$. However, recloser-recloser's coordination margin is interfered, $(T_{R,j}-T_{R,j+1})^2 <{\Delta} T_{RR}$}
\label{fig: adplimit}
\vspace{-4mm}
\end{figure}

Existing solutions can address the protection problem under limited DG penetration. However, at high DG penetration level, solely adjusting recloser/relays' settings will not be sufficient to retrieve protection coordination. An example is illustrated by Fig.~\ref{fig: adplimit}. The fuse-recloser coordination is disrupted due to a large fault current contribution from DG. Shifting down the TCI curve of Recloser $j$ will restore the coordination margin between the fuse and recloser, but it will consequently interferes with the coordination margin between Recloser $j$ and Recloser $j+1$. This problem cannot be solved by adjusting other TCI curves, because the fuse's TCI curve is inflexible to be adjusted, and the downstream Recloser $j+1$ might have reached its minimum time setting.   

\section{Modeling DG's fault current contribution}
Due to the limitation of recloser/relay's adaptivity, we consider retrieving protection coordination through controlling DG's fault current contribution in active distribution systems. DG's fault current control can be achieved through hardware deployment, such as fault current limiter \cite{tang2005application, El-Khattam}. However, this method requires extra capital investment and lacks of implementation flexibility. An alternative way is to control DG's output power in pre-fault networks, whose principle is reasoned in this section. 
\vspace{-2mm}

\subsection{DG's types} 

\begin{figure}[ht!]
\setlength{\abovecaptionskip}{-3 pt}
\centering
\subfigure[Synchronous]{
\hspace{-4mm}
  \includegraphics[width= .2\textwidth]{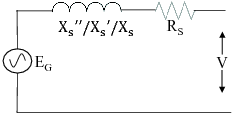}
  \label{fig:1a}}
\subfigure[Inverter-based]{
  \includegraphics[width= .09\textwidth]{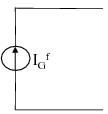}
  \label{fig:8a}}
 \subfigure[Asynchrounous]{
  \includegraphics[width= .35\textwidth]{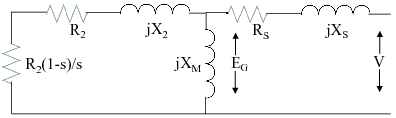}
  \label{fig:1b}}\\
\caption{The\'venin equivalent circuit for three major DG types}

\label{fig: DGtypes}
\end{figure}

There are three major DG types: synchronous generator, asynchronous generator, and inverter-based generator. Their circuits are shown in Fig.~\ref{fig: DGtypes}. 

Synchronous DG includes those driven by conventional energy sources and by renewable energies, such as biomass and solar steam systems \cite{nimpitiwan2007fault}. Their fault current contribution depends on the generator's parameters (subtransient and transient reactance) as well as the equivalent pre-fault voltage. 

Asynchronous generators are represented by Type I and Type II wind turbine generators (i.e. squirrel cage generator and wound rotor induction generator) \cite{windpsrc}. Their initial fault response is determined by the pre-fault slip, which is determined by the wind speed and power output before fault happens. This fault current lasts for a few cycles and diminishes.  

Most DG driven by renewable energy sources are inverter-based, such as, Type III and Type IV wind turbine generators (i.e. variable speed double fed generator and induction generator interfaced through a full AC/DC/AC converter), Photo-Voltaic (PV) panels, fuel cells and etc. \cite{hossain2014, Baran}. The inverter control can be current-based or voltage-based. Current-based control is more often used\cite{Plet,haron2012review}. During the fault response, DG unit's peak fault current occurs during the subtransient timeframe and is determined by its pre-fault terminal voltage. If this subtransient fault current is greater than a  threshold ($2-3$ times of a DG unit's rating current), the inverter is turned off and no fault current is flown from DG. After several cycles, the interter clamps the fault current  to a fix value, which is $1.25-2$ times of a DG unit's rating current \cite{Plet, Baran}.
\vspace{-2mm}
\subsection{Relationship between DG's fault current and power output}
In active distribution systems, DG's pre-fault voltage is determined by feeder's load flow. According to the original version of IEEE 1547, DG shall not actively regulate the voltage at its connection point \cite{ieee1547}. In the 2014 amendment, DG is allowed to regulate voltage under approval of the system operator by changes of its real and reactive power \cite{basso2004ieee}. Based on either standard, DG's terminal voltage is dependent of its output power. Therefore, in calculation of network's pre-fault voltages, DG should be modeled as PQ bus rather than PV bus. This is the main difference between the fault analysis in transmission and active distribution systems. DG's pre-fault terminal voltage can be found by solving the load flow on feeders from\cite{ffwu},
\begin{align}\label{distflow}
\begin{split}
P_{i+1}= P_{i} -r_i(P^2_i+Q^2_i)/V^2_i - (P_{D,i} -P_{G,i})\\ 
Q_{i+1}= Q_{i} -x_i(P^2_i+Q^2_i)/V^2_i - (Q_{D,i} -Q_{G,i})
\end{split}
\end{align}  

After $V_i$ is obtained, DG's Th{\'e}venin equivalent voltage and impedance, for synchronous and asynchronous generators, can be calculated from the circuits in Fig.~\ref{fig: DGtypes}. Inverter-based DG is treated as a constant current source. Then DGs' fault current contribution is solved the combined Th{\'e}venin equivalent circuits of the whole distribution network.     
Directional recloser/relays only pick up fault currents from their upstream. Given that DG is only connected to feeders, lateral fuses are downstream to all DG units. During a fault, a fuse sees the fault current from all DG units and the substation. The recloser to be coordinated with the fuse at node $j$ only sees fault current from the substation and DG upstream. The fault current disparity between the fuse and recloser is contributed from all DG units downstream node $j$, expressed as   
\begin{equation} \label{eq:recfusedisp}
\Delta I_{FR,j} = \sum_{i\in\mathcal{I}^u_j} I^f_{G,i} .
\end{equation}  

Fault current disparity between Recloser $j$ and Recloser $j+1$ is determined by the fault current contribution from the DG on the section between node $j$ and $j+1$, expressed as:     
\begin{equation} \label{eq:recrecdisp}
I_{RR,j}^\Delta = \sum_{i \in\mathcal{I}^j} I_{G,i}^f .
\end{equation}

%
%
%

%

%
%

\section{Proposed Protection Scheme}
In this section, we propose a new protection method that simultaneously adjusts recloser/relays' setting and DG's output power in pre-fault networks. The proposed method contains two objectives, minimum delay of backup protective devices and maximum DGs' total output power. The challenge of solving the optimization problem arises from the high non-linearity of recloser/relays' TCI inverse curves. An algorithm is proposed to address this challenge.  

\vspace{-2mm}
\subsection{Formulation}
The first objective is to minimize backup delay for each pair of reclosers. For reclosers on the same feeder, theri TCI curves have similar shapes. Therefore, the reclosers' coordination can be visualized by stacking the TCI curves of all recloser/relays from the one at the feeder's end to that at the feeder's head. Therefore, a minimum total clearing time must lead to minimum backup delay for each pair \cite{anderson}~\cite{Abdelaziz}. The objective is equivalent to, 

\begin{equation} \label{eq:mincleartime}
\min \sum_{\mathcal J} T_{R,j}
\end{equation} 

The second objective is to maximize DG's output, 

\begin{equation} \label{eq:maxDG}
\max \sum_{\mathcal{I}} P_{G,i} ,
\end{equation}


The two objectives, \eqref{eq:mincleartime} and \eqref{eq:maxDG} are bounded by the protection coordination conditions \eqref{cod_range} and \eqref{cod_margin}. For fuse-recloser and recloser-recloser coordinations, these conditions are specified as, 
\begin{equation} \label{eq:fusereccoor}
T_{F,j} - T_{R,j} \geq \Delta T_{FR},
\end{equation}
\begin{equation} \label{eq:recreccoor}
T_{R,(j-1)} - T_{R,j} \geq \Delta T_{RR}.
\end{equation}


Since the reclosing sequence could contain more than one fast tripping curve, \eqref{eq:fusereccoor} and \eqref{eq:recreccoor} are applied according to the anti-islanding requirements and DG's protection settings. If DG is disconnected from the grid after the recloser's first trip, then \eqref{eq:mincleartime} only coordinates recloser's first fast tripping curve and fuse's MM curve, and \eqref{eq:maxDG} coordinates two reclosers' first fast tripping curves.  

The TCI curve for recloser/relays is specified by \cite{anderson}, 
\begin{equation}\label{eq:TCchar}
T_{R,j} = \frac{a\cdot D_{R,j}}{\left(\frac{I_{R,j}^f}{I^p_{R,j}}\right)^m - c}+b\cdot D_{R,j}+K ,
\end{equation} 
where  $a,b,c,m,$ and $K$ are constants provided by manufactures. 




\vspace{-2mm}
\subsection{Algorithm}
Equation \eqref{eq:TCchar} is highly nonlinear of the recloser's fault current $I_{R,j}^f$. To reduce the nonlinearity, the formulation is decomposed into two sub-problems, whose objectives are minimizing total fault clearing time and maximizing DG's output. The first sub-problem is constrained by \eqref{eq:fusereccoor} and \eqref{eq:recreccoor}. The second sub-problem is constrained by \eqref{eq:fusereccoor}. To further simplify the second sub-problem, constraint \eqref{eq:fusereccoor} is transformed into a function of DGs' fault current contribution as, 
\begin{equation}\label{eq:fuserecmargin}
\Delta I^f_{FR,j} \leq I^f_{max, F,j} - I^f_{max, R,j} -\Delta I_{FR}
\end{equation}

where $\Delta I_{FR}$ is the current value corresponding to the fuse-recloser coordination margin $\Delta T_{FR}$. Left hand side of \eqref{eq:fuserecmargin} can be obtained from \eqref{eq:recfusedisp}. The maximum fault currents and  $\Delta I_{FR}$ on right hand side are obtained from the first sub-problem. 

%

\begin{figure}[!t]
\centering
\includegraphics[width=3in]{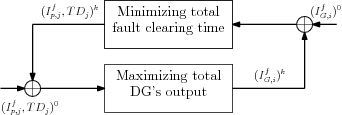}
\caption{Algorithm for the proposed protection scheme. The decomposed optimization problems are solved alternatively by updating the coupling constraint. The initial value is input as $(.)^0$, the kth iteration is denoted as $(.)^k$.}
\label{fig: Fig6}
\vspace{-4mm}
\end{figure}

Based on the decomposed mathematical structure, we propose an algorithm, illustrated in Fig.~\ref{fig: Fig6}. The two objectives are optimized alternatively. Each sub-problem is optimized after updating the coupling constraints, \eqref{eq:fusereccoor} and \eqref{eq:fuserecmargin}, with the optima obtained from the other problem. By doing so, the non-linearity are extinguished in the first sub-problem, greatly reducing the computation size and complexity.
%

The algorithm may start with either sub-problem. For example, the solver starts with the second subproblem, maximizing DGs' output. Initial settings are picked for TCI curves, which implies the maximum fault current in the system. This set of initial values, $(I^p_{R,j},D_{R,j}, I^f_{max, F,j},I^f_{max, R,j})^0$, may be chosen heuristically or from undergoing system operation conditions. The solver calculates for DG's optimal output power $(P_G,i)^1$. With this value, reclosers' fault current $(I^f_{R,j})^1$ can be found by following the procedure in Section III. Then constraint \eqref{eq:fusereccoor} is updated in the first sub-problem. The solver optimizes the total clearing time is minimized and obtains optimal recloser/relay settings, $(I^p_{R,j},D_{R,j})^1$. This completes the first iteration. The iterative optimization process ends when the constraints in the sub-problems cannot be relaxed further by controlling DG's output and adjusting recloser settings.

\subsection{Implementation and Dependence on Distribution Automation} 
The proposed protection method and algorithm enable implementation flexibility based on system operation conditions. In practice, the frequencies of adapting reclosers and controlling DG depend on a few factors, such as DGs' type and intermittency, authorized intervention to DG's natural output, communication readiness of the grid, and etc. For DG of high intermittency and their fault current highly dependent of power output, such as Type I and II wind turbine generator and solar steam systems, the protection settings could be adapted hourly or more frequently. This flexibility is tested in Section V with practical numerical examples. 

The proposed protection method requires metering at all DG locations and remote control of the authorized DG units. Today's distribution grid may not completely fulfill these requirements. The proposed method provides the flexibility to be implemented progressively in accordance with the grid's available automation level. Equation \eqref{eq:recfusedisp} to \eqref{eq:recrecdisp} show that the requirements of metering and communication declines as the protection goals become more critical. For example, for fuses' secure operation, only DG at reclosers' upstream need to be controlled; for minimum backup delay, DG need to be controlled at every feeder section between two reclosers. The later protection goal requires more discrete monitoring and communication, but presents relatively less protection criticalness. This is a favorable characteristic in implementation, in the sense that the most critical protection requirement can be fulfilled at an earlier stage of distribution automation.

\section{Numerical Examples}
\begin{figure}[!t]
\vspace{-3mm}
\centering
\includegraphics[width=0.25\textwidth,angle=90]{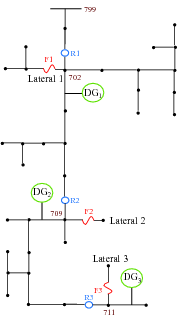}
\caption{IEEE 37-bus distribution system with DGs connected in the main feeder.}
\label{fig: system}
\vspace{-4mm}
\end{figure}

This section demonstrates the benefits of the proposed protection scheme with an IEEE 37-bus system. The system configuration is altered with DG integration and shown in Fig.~\ref{fig: system}.  The main feeder is protected with 3 reclosers, which locations are shown in Table~\ref{table:faultcurr}. The proposed protection method, short circuit analysis and power flow were implemented in MATLAB.

In the test system, the reclosers' TCI characteristics are modeled with \eqref{eq:TCchar}. Without losing generality, the reclosers' Time-Dial Setting (TDS) and pickup current are assumed to be continuous. The TDS' are selected from the interval, ${D}_{R,j} \in [0.1,1]$, and the pickup currents $I^p_{R,j}$ are chosen twice greater than the maximum load current and less than 50\% of Line-Line fault current~\cite{anderson}. Fuses' TC characteristics are represented by their minimum melting curves and are linearized ~\cite{anderson}.

The initial settings of the protective devices are based on the configuration without DG's integration. The maximum currents at reclosers' bus are calculated in Table~\ref{table:faultcurr}. The reclosers are set based on the procedure described in Section II. Since R3 only needs to coordinate with its lateral fuse, its TDS can be selected to be the lowest, that is ${D}_{R,3}=0.1$. The rest of the reclosers are coordinated in pairs.  

\begin{table}[ht]
\vspace{-5mm}
\caption{Reclosers' Max/Min Fault Currents}
\centering
\begin{tabular}{|c |c | c |}
\hline 
Recloser & Bus & Max  Fault Current (A) \\ [0.5ex]
\hline
R1 & 702 & 2975   \\ [0.5ex]
\hline
R2 & 709 & 2445   \\ [0.5ex]
\hline
R3 & 711 & 1823   \\ [0.5ex]
\hline
\end{tabular}
\label{table:faultcurr}
\end{table}
\vspace{-5mm}

\subsection*{Case A: Single-step implementation}
\begin{figure}[!t]
\setlength{\abovecaptionskip}{-3 pt}
\centering
\includegraphics[width=0.5\textwidth]{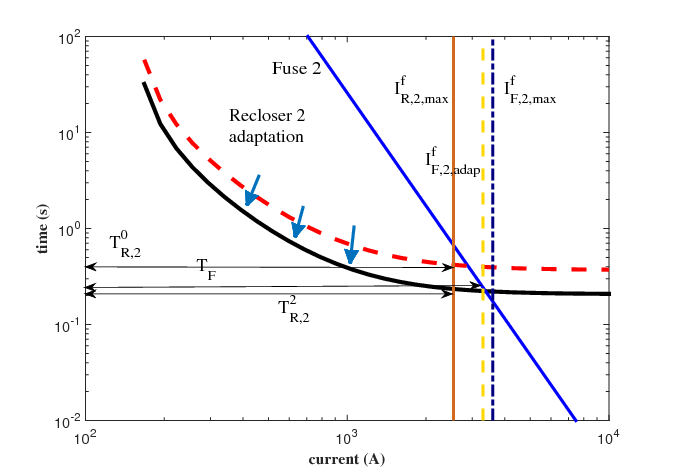}
\caption{TCI curves of R2 and F2. R2 Original and R2 adaptive are the R2's TC curves using basic (non-adaptive) protection and adaptive protection; initial and adapted DG power output.}
\label{fig: Fig14r}
\vspace{-4mm}

\end{figure}

The proposed method is tested assuming all DG are controllable and their output known in advance. For a 3-phase fault on Lateral 2, DG's integration will disrupt coordination between the recloser, R2, and fuse, F2, as shown in Fig.~\ref{fig: Fig14r}. Due to the large high penetration level of DG in the systme, coordination cannot be restored by solely adapting R2's TCI curve. With the proposed method, DG's output is controlled to regain the coordination margin. The tripping sequence restoration procedure is illustrated in Fig.~\ref{fig: Fig13r}.


\begin{figure}[!t]
\setlength{\abovecaptionskip}{-3 pt}
\centering
\includegraphics[width=0.5\textwidth]{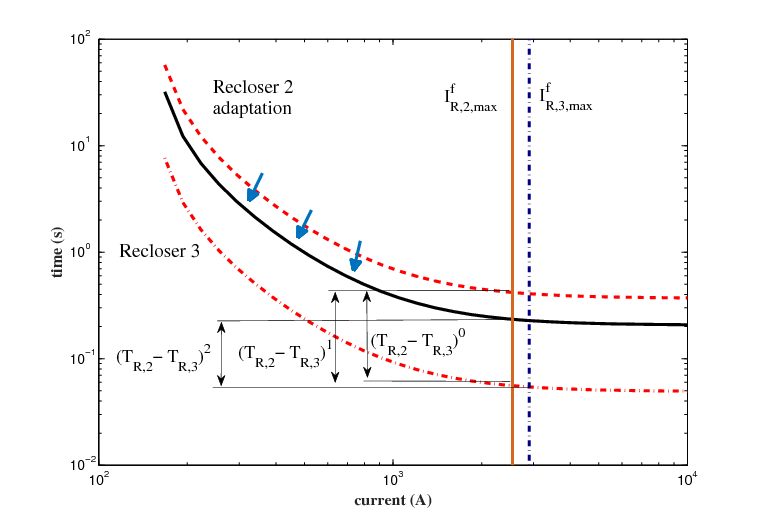}
\caption{TCI curves of R2 and R3.R2 Original and R2 adaptive are the R2's TC curves using basic (non-adaptive) protection and adaptive protection, respectively. }
\label{fig: Fig12r}
\vspace{-4mm}
\end{figure}

\begin{figure}[!t]
\setlength{\abovecaptionskip}{-3 pt}
\centering
\includegraphics[width=0.5\textwidth, height=7cm]{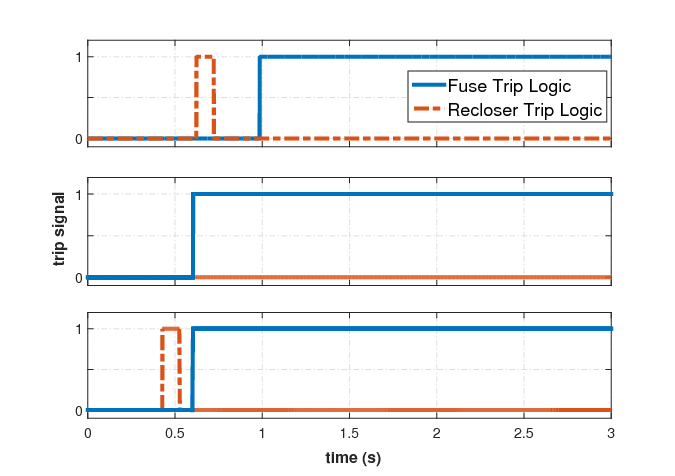}
\caption{Tripping sequence of F2 and R2 in Fig.~\ref{fig: system}. Top: no DG is connected in the test system; Middle: with DG integrated and basic (non-adaptive) protection; Bottom: with DG integrated and conventional adaptive protection}
\label{fig: Fig13r}
\vspace{-4mm}
\end{figure}

Reclosers' coordination is examined under basic and adaptive protection. Recloser R2 backups R3 for a fault on Lateral 3. Fig.~\ref{fig: Fig12r} shows that, before DG's integration, the backup time interval between R2 and R3 is $T_0 = (T_{R,2}-T_{R,3})^0$. DG2's connection increase the maximum fault current seen by R3 as $I^{f}_{R,3,\max}$. The resultant delay between R2 and R3 is increased to $T_1 = (T_{R,2}-T_{R,3})^1$ and $T_0 < T_1$. To reduce this backup clearing time, R2 is adapted to reduce the backup fault clearing time to $T_2 = (T_{R,2}-T_{R,3})^2$.


\subsection*{Case B: Multiple-step implementation }
\begin{figure}[!t]
\setlength{\abovecaptionskip}{-3 pt}
\centering
\includegraphics[width=0.5\textwidth]{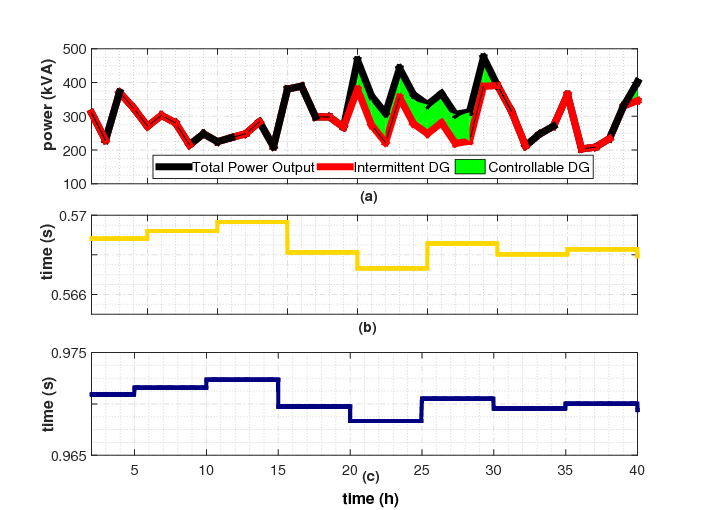}
\caption{Proposed protection scheme under DG of mixed types integration. The test system is shown in Fig.~\ref{fig: system}. (a) DG's output; (b) Time-Dial Setting of R2; (c) total fault clearing time of all reclosers.}
\label{fig: Fig16r}
\vspace{-4mm}
\end{figure}

The benefits of proposed protection method is fully realized when DG of mixed types is integrated to the system. As shown in Fig.~\ref{fig: Fig16r}(a), both inverter-based and controllable DG present in the test system. We demonstrate the flexibility of the proposed scheme by adapting DG's output and recloser settings at one hour and five hours, respectively. This scenario has practical implications, since protective systems are less often adapted and DG, as could be involved in demand response programs, is not rare controlled hourly. 

Fig.~\ref{fig: Fig16r} shows that the total DG output is maximized to reduce the fault clearing time in the system. At peak period of renewable generation, the controllable DG units are commanded to produce less power, to ensure recloser and fuse coordination and prevent unnecessary fuse operation, as in \eqref{eq:fusereccoor} and \eqref{eq:recfusedisp}.   

\section{Conclusion}
Legacy protection system at distribution level assumes power flowing in single direction from the substation. Increasing penetration of Distributed Generation (DG) invalidates this assumption and will disrupt protection coordination. Solutions proposed in previous studies include (i) limiting DG's fault current contribution; (ii) deploying microprocessor-based recloser/relays; and (iii) deploying non-overcurrent protection schemes. Deficiencies of these solutions can be attributed to their high cost, inflexibility of implementation, and limited effectiveness under high penetration of renewable DG.

To retrieve protection coordination in active distribution grids, this paper proposes a protection method which integrally adjusts the settings of protection devices and DG's output in pre-fault distribution networks. The formulation of the proposed method is decomposed into two optimization problems, aiming at two distinct objectives: maximizing total power output from DG's, and retaining protection coordination and sensitivity. This decomposed mathematical structure embodies the following benefits: (i) reducing computation effort by extinguishing the non-linearity of relay/reclosers' time-current inverse characteristics; (ii) enabling adaptive protection and DG control under different timeframes; (iii) providing implementation flexibility based on available communication and automation level of distribution grids. These benefits are demonstrated through numerical cases on the IEEE-37 bus system.

\bibliographystyle{ieeetr}
\bibliography{references}

\begin{thebibliography}{10}

\bibitem{Girgis}
A.~Girgis and S.~Brahma, ``Effect of distributed generation on protective
  device coordination in distribution system,'' in {\em Power Engineering,
  2001. LESCOPE '01. 2001 Large Engineering Systems Conference on},
  pp.~115--119, 2001.

\bibitem{Brahma}
S.~Brahma and A.~Girgis, ``Microprocessor-based reclosing to coordinate fuse
  and recloser in a system with high penetration of distributed generation,''
  in {\em Power Engineering Society Winter Meeting, 2002. IEEE}, vol.~1,
  pp.~453--458 vol.1, 2002.

\bibitem{Chaitusaney}
S.~Chaitusaney and A.~Yokoyama, ``Prevention of reliability degradation from
  recloser-fuse miscoordination due to distributed generation,'' {\em Power
  Delivery, IEEE Transactions on}, vol.~23, pp.~2545--2554, Oct 2008.

\bibitem{ieee2004impact}
``Impact of distributed resources on distribution relay protection {IEEE}
  {P}ower {S}ystem {R}elay {C}ommittee,'' Aug 2004.

\bibitem{de2012impact}
S.~De~Bruyn, J.~Fadiran, S.~Chowdhury, S.~Chowdhury, and P.~Kolhe, ``The impact
  of wind power penetration on recloser operation in distribution networks,''
  in {\em Universities Power Engineering Conference (UPEC), 2012 47th
  International}, pp.~1--6, IEEE, 2012.

\bibitem{tang2005application}
G.~Tang and M.~Iravani, ``Application of a fault current limiter to minimize
  distributed generation impact on coordinated relay protection,'' in {\em Int.
  Conf. Power Systems Transients, Montreal, QC, Canada}, pp.~19--23, 2005.

\bibitem{El-Khattam}
W.~El-Khattam and T.~Sidhu, ``Restoration of directional overcurrent relay
  coordination in distributed generation systems utilizing fault current
  limiter,'' {\em Power Delivery, IEEE Transactions on}, vol.~23, pp.~576--585,
  April 2008.

\bibitem{Abdelaziz}
A.~Abdelaziz, H.~Talaat, A.~Nosseir, and A.~A. Hajjar, ``An adaptive protection
  scheme for optimal coordination of overcurrent relays,'' {\em Electric Power
  Systems Research}, vol.~61, no.~1, pp.~1 -- 9, 2002.

\bibitem{dewadasa2011protection}
M.~Dewadasa, A.~Ghosh, and G.~Ledwich, ``Protection of microgrids using
  differential relays,'' in {\em Universities Power Engineering Conference
  (AUPEC), 2011 21st Australasian}, pp.~1--6, IEEE, 2011.

\bibitem{haron2012review}
A.~R. Haron, A.~Mohamed, and H.~Shareef, ``A review on protection schemes and
  coordination techniques in microgrid system,'' {\em Journal of Applied
  Sciences}, vol.~12, no.~2, p.~101, 2012.

\bibitem{sortomme2010microgrid}
E.~Sortomme, S.~Venkata, and J.~Mitra, ``Microgrid protection using
  communication-assisted digital relays,'' {\em Power delivery, IEEE
  transactions on}, vol.~25, no.~4, pp.~2789--2796, 2010.

\bibitem{Barker}
P.~Barker and R.~de~Mello, ``Determining the impact of distributed generation
  on power systems. i. radial distribution systems,'' in {\em Power Engineering
  Society Summer Meeting, 2000. IEEE}, vol.~3, pp.~1645--1656 vol. 3, 2000.

\bibitem{chaitusaney2005impact}
S.~Chaitusaney and A.~Yokoyama, ``Impact of protection coordination on sizes of
  several distributed generation sources,'' in {\em Power Engineering
  Conference, 2005. IPEC 2005. The 7th International}, pp.~669--674, IEEE,
  2005.

\bibitem{Fumilayo}
H.~B. Funmilayo and K.~L. Butler-Purry, ``An approach to mitigate the impact of
  distributed generation on the overcurrent protection scheme for radial
  feeders,'' in {\em Power Systems Conference and Exposition, 2009. PSCE '09.
  IEEE/PES}, pp.~1--11, March 2009.

\bibitem{anderson}
P.~Anderson, {\em Power System Protection}.
\newblock IEEE Press Series on Power Engineering, Wiley, 1998.

\bibitem{Kersting}
W.~Kersting, ``Radial distribution test feeders,'' {\em Power Systems, IEEE
  Transactions on}, vol.~6, pp.~975--985, Aug 1991.

\bibitem{glover2011power}
J.~D. Glover, M.~Sarma, and T.~Overbye, {\em Power System Analysis \& Design,
  SI Version}.
\newblock Cengage Learning, 2011.

\bibitem{Fazanehrafat}
A.~Fazanehrafat, S.~Javadian, S.~Bathaee, and M.-R. Haghifam, ``Maintaining the
  recloser-fuse coordination in distribution systems in presence of dg by
  determining dg's size,'' in {\em Developments in Power System Protection,
  2008. DPSP 2008. IET 9th International Conference on}, pp.~132--137, March
  2008.

\bibitem{basso2004ieee}
T.~S. Basso and R.~DeBlasio, ``Ieee 1547 series of standards: interconnection
  issues,'' {\em Power Electronics, IEEE Transactions on}, vol.~19, no.~5,
  pp.~1159--1162, 2004.

\bibitem{kumpulainen2004analysis}
L.~Kumpulainen and K.~Kauhaniemi, ``Analysis of the impact of distributed
  generation on automatic reclosing,'' in {\em Power Systems Conference and
  Exposition, 2004. IEEE PES}, pp.~603--608, IEEE, 2004.

\bibitem{nimpitiwan2007fault}
N.~Nimpitiwan, G.~T. Heydt, R.~Ayyanar, and S.~Suryanarayanan, ``Fault current
  contribution from synchronous machine and inverter based distributed
  generators,'' {\em Power Delivery, IEEE Transactions on}, vol.~22, no.~1,
  pp.~634--641, 2007.

\bibitem{windpsrc}
``Fault current contributions from wind plants,'' in {\em Protective Relay
  Engineers, 2015 68th Annual Conference for}, pp.~137--227, March 2015.

\bibitem{hossain2014}
J.~Hossain and A.~Mahmud, {\em Large Scale Renewable Power Generation: Advances
  in Technologies for Generation, Transmission and Storage}.
\newblock Green Energy and Technology, Springer, 2014.

\bibitem{Baran}
M.~E. Baran and I.~El-Markaby, ``Fault analysis on distribution feeders with
  distributed generators,'' {\em IEEE Transactions on Power Systems}, vol.~20,
  pp.~1757--1764, Nov 2005.

\bibitem{Plet}
C.~A. Plet, M.~Graovac, T.~C. Green, and R.~Iravani, ``Fault response of
  grid-connected inverter dominated networks,'' in {\em IEEE PES General
  Meeting}, pp.~1--8, July 2010.

\bibitem{ieee1547}
``Ieee application guide for ieee std 1547(tm), ieee standard for
  interconnecting distributed resources with electric power systems,'' {\em
  IEEE Std 1547.2-2008}, pp.~1--217, April 2009.

\bibitem{ffwu}
M.~Baran and F.~F. Wu, ``Optimal sizing of capacitors placed on a radial
  distribution system,'' {\em IEEE Transactions on Power Delivery}, vol.~4,
  pp.~735--743, Jan 1989.

\end{thebibliography}

\begin{IEEEbiography}{J.K. Wang}
Jiankang (J.K.) Wang (M’2014, SM’2007) received her M.S. and PH.D. degree from the Department of Electrical Engineering and Computer Science at Massachusetts Institute of Technology, in 2009 and 2013, respectively. Her research focus is on power distribution systems’ operation and planning, end-user engagement through electricity markets. She is now an assistant professor in the Department of Electrical and Computer Engineering, at the Ohio State University. She has a joint appointment from the Department of Integrated System Engineering.
\end{IEEEbiography}
\begin{IEEEbiography}{Christian Moya Calderon}
(S'14) received the B.Sc. degree in electronic engineering from Universidad San Francisco de Quito, Ecuador, in 2010 and the M.Sc. degree in electrical and computer engineering from The Ohio State University, OH, USA, in 2014.
He is currently a Ph.D. student in the department of electrical and computer engineering at The Ohio State University, OH, USA. His main research interest include power networks and control systems. 
\end{IEEEbiography}
\end{document}